\documentclass[english,a4paper,11pt]{article}
\usepackage{theorem}
\usepackage{amssymb,amsmath}
\usepackage{latexsym}
\usepackage{babel}
\usepackage{color}
\usepackage{slashbox}
\usepackage[usenames,dvipsnames,svgnames,table]{xcolor}
\usepackage{color,graphicx,array,fleqn}
\usepackage{bbm}
\usepackage{subfig}

\newtheorem{theorem}{Theorem}[section]

\newtheorem{corollary}[theorem]{Corollary}
\newtheorem{lemma}[theorem]{Lemma}
\theorembodyfont{\normalfont\mdseries\upshape}

\newtheorem{example}[theorem]{Example}
\newtheorem{definition}[theorem]{Definition}

\title{\bf Numerical Solution for an Inverse Variational Problem}
\author{A.I. Garralda-Guillem$^{\hbox{(a)}}$ and P. Montiel L\'opez$^{\hbox{(b)}}$}

\begin{document}

\maketitle

\centerline{$^{\hbox{(a)}}$University of Granada,
E.T.S. Ingenier\'{\i}a de Edificaci\'on,}

\centerline{ Department of Applied Mathematics, Granada (Spain), e-mail:
agarral@ugr.es}

\medskip

\centerline{$^{\hbox{(b)}}$University of Granada,
Centro de Magisterio La Inmaculada,} 

\centerline{Department of Sciences, Granada (Spain), e-mail: pablomontiel@cmli.es}

\begin{abstract}
\noindent In the present work, firstly, we use a minimax equality to prove the existence of a solution of certain system of varitional equations and we provide a numerical approximation of such a solution. Then, we propose a numerical method to solve a collage-type inverse problem associated with the corresponding system, and illustrate the behaviour of the method with a numerical example. 
\end{abstract}



\maketitle

\section{Introduction}

\noindent In this paper, we propose a numerical method to solve an inverse problem associated with a certain system of variational equations, which is based upon a generalization of the classical collage theorem. To this end, we first characterize, in terms of the existence of a scalar, the solvability of the following system
\[
\
\left\{\begin{array}{lll}

x^{*}_1(y_1-x_0) & \leq & a_1(y_1,y_1-x_0)\\
x^{*}_2(y_2-x_0) & \leq & a_2(y_2,y_2-x_0)\\

\; \; \; \; \; \; \; \; \; \;  \vdots & \;  \vdots & \; \; \; \; \; \; \; \; \; \; \vdots\\

x^{*}_N(y_N-x_0) & \leq & a_N(y_N,y_N-x_0)\\

\end{array},
\right.
\]
with $x^{\ast}_i$ being continuous and linear functionals in a real reflexive Banach space and $a_i$ continuous bilinear forms in the same space.

Since Stampacchia's results in the 1960s, the study of variational inequalities and systems of variational equations has raised great interest, in part due to the fact that a wide range of optimization problems can be reformulated as variational problems. The concept of variational systems encompasses different types of problems, for example, in \cite{pan}, the Nash equilibrium problem, the spatial equilibrium problem and the general equilibrium programming problem are modeled as a system of variational inequalities. In \cite{garrui2019}, the variational system includes certain mixed variational formulations associated with some elliptic problems. Also, we can find them associated with abstract economy, \cite{ans99}.

Another specific type problem of variational systems, and that is related to our work, is the so-called common solutions to variational inequalities problem, which consists of finding common solutions to a system of variational inequalities. There are different approaches to this kind of inequality systems: in \cite{zao10} the definition domains of the functions of the system are closed and convex sets of a Hausdorff topological vector space. In \cite{kas00}, the system problem is dealt with for only two inequalities, a treatment that is generalized in \cite{cen12}.

In the study of the solution of variational equations or systems of variational equations, a wide range of techniques is used, including those of minimax type, \cite{BiCasPa, Fan1972, Par2018}, or those that use fixed point results and their associated iterative methods, such as those we detail next. In \cite{ans99}, the authors prove the existence of solution of certain variational systems by using a multivalued fixed point theorem. Also, one can find a proof of the existence of solution to a variational system with the Brouwer fixed point theorem in \cite{zao10} as well as the construction of an iterative algorithm for approximating the unique solution of the system and a discussion of the convergence analysis of the algorithm.

In the present paper, we use a minimax technique to prove the existence of solution for the system of variational inequalities, Theorem \ref{th:1}, that, unlike the different results that appear in the articles mentioned above, characterizes the existence of solutions in closed subsets (not necessarily convex).

Once the conditions that ensure the existence of a solution for the system of variational equations have been established, we will deal with the inverse problem, i.e., assuming that the model which depends on different parameters has been established, and some empirical solutions have been obtained, we will try to approximate the parameters for which the empirical solutions obtained are an approximation of the solution of the theoretical model.

From the different approaches proposed in the literature to solve inverse problems, we rely on the approach of the so-called Collage theorem, which starts by considering the forward problem as a solution to a fixed point problem and deduces its analysis from Banach's fixed point theorem. We follow the line of different proposed generalizations of the Collage theorem, which are supported by different versions of the Lax-Milgram theorem, established, for example, to solve inverse problems associated with different families of integral or ordinary differential equations, (\cite{Capasso2014, KuVr99, kunze03a, kunze03b}), or of partial differential equations (\cite{berkunlatrui, garrui2019, garrui2014, KuLaToVr06-sub, Ku1, KuLaToLeRuGa15, Levere}).

The paper is organized as follows. The first section begins with the presentation of our minimax tool and the variational system. Theorem \ref{th:1} is the central point of this section, and provide us with a characterization of the solvability of the variational system. Moreover, from this theorem we derive a result which implies Stampacchia's theorem. The following section begins with a collage-type result that will be used in the numerical treatment of the inverse problem of a concrete example. To this end, we first propose a numerical approximation of the solution of the forward problem that we show in different tables and graphics. Finally, we finish our work with the conclusions.

\section{The forward variational problem}

In this section we deal with a result, Theorem \ref{th:1}, that generalizes the classic Stampacchia theorem. Indeed, it allows us to characterize the existence of a solution to a system of variational equations as that of a certain scalar. We should mention that minimax inequalities are a widely used technique in variational analysis: \cite{aub1998} is a good example. In \cite{sim1998}, we see the equivalence between minimax results and the Hahn-Banach Theorem, and how these results are used as functional analytic tools.

The fundamental tool to establish this direct result is given by the following minimax inequality \cite{kaskol1996}, which includes a not very restrictive convexity condition, which allows to characterize the validity of the minimax identity \cite{rui2014, rui2016}. This concept of weak convexity is called infsup-convexity, and it appears with a  nomenclature for first time as affine weakly convexlikeness in \cite{ste}. The infsup-convexity arises, in a natural way, when we deal with equilibrium and minimax problems.

\begin{definition}
If $X$ and $Y$ are nonempty sets, a function
$g: X \times Y  \longrightarrow \mathbb{R}$ is called
\textit{infsup-convex} on $Y$ provided that
\[
\inf_{y \in Y} \max_{x \in X} g(x,y) \le \max_{x \in X}
\sum_{j=1}^m t_j g(x, y_j),
\]
whenever $m \ge 1$, $y_1,\dots,y_m \in Y$ and $t \in \Delta_m$, where $\Delta_m$ is the probability simplex, $\Delta_m:= \{ (t_1,\dots,t_m) \in \mathbb{R}^m : \ t_1,\dots,t_m \ge 0  \hbox{ and } \displaystyle \sum_{j=1}^m t_j =1 \}$.
\end{definition}

Clearly, the infsup-convexity extends the concept of convex function, but also another types of weak convexity, such as convexlikeness (\cite{fan}). Let us recall that a function $f: X \times Y \longrightarrow \mathbb{R}$ in convexlike (or Fan-convex) on $Y$ when for any $y_1,y_2\in Y$ there exist $y \in Y$ and $0<t<1$ such that 
\[
x \in X \ \Rightarrow \ f(x,y) \le tf(x,y_1)+(1-t)f(x,y_2).
\]

The concept of infsup-convexity is used in the following minimax result, \cite{kaskol1996}.

\begin{theorem}\label{th:minimax}
Assume that $X$ is a nonempty, convex and compact subset of a real topological vector space, $Y$ is a nonempty set and $g: X \times Y \longrightarrow \mathbb{R}$ is continuous and concave on $X$. Then, 
\[
\inf_{y \in Y} \max_{x \in X} g(x,y) = \max_{x \in X} \inf_{y \in Y} g(x,y)
\]   
if, and only if, $g$ is infsup-convex on $Y$.

\end{theorem}

\bigskip

This minimax inequality is of the Hahn-Banach type, in the sense that it is equivalent to this central result of the functional analysis. In fact, the Hahn-Banach theorem and some of its generalizations have also been used to prove some variational results \cite{sai2018, sim2007}, even for some systems of variational equations \cite{garrui2019, garrui2014} that include, as a particular case, those corresponding to the mixed variational formulation of the classical Babu\v{s}ka-Brezzi theory \cite{bofbrefor2013, gat2014}.

Now, we introduce in a precise way the forward problem involving a system of variational equations. Let $E$ be a real and reflexive Banach space, let  $Y_{1}, \dots, Y_{N}$ be closed and nonempty subsets of $E$, $\displaystyle Y :=\prod_{i=1}^N Y_i $, let $x^{*}_1:E \longrightarrow \mathbb{R},  \dots, \, x^{*}_N:E \longrightarrow \mathbb{R}$ be continuous and linear functionals and let $a_1:E \times E \longrightarrow \mathbb{R}, \dots, \, a_N:E \times E \longrightarrow \mathbb{R}$ be continuous bilinear forms. We consider the following variational problem: find an $\displaystyle x_0 \in \overline{Y}:= \bigcap^N_{i=1} Y_i$ such that 
\begin{equation}\label{eq:1}
y \in Y \ \Rightarrow \ 
\
\left\{\begin{array}{lll}

x^{*}_1(y_1-x_0) & \leq & a_1(y_1,y_1-x_0)\\
x^{*}_2(y_2-x_0) & \leq & a_2(y_2,y_2-x_0)\\

\; \; \; \; \; \; \; \; \; \;  \vdots & \;  \vdots & \; \; \; \; \; \; \; \; \; \; \vdots\\

x^{*}_N(y_N-x_0) & \leq & a_N(y_N,y_N-x_0)\\

\end{array}.
\right.
\end{equation}

In order to study this system, we note by  $x^{*}$ the linear and continuous functional defined in $ Y$ as
\[
 x^{*}(y):= x^{*}_1(y_1)+ \cdots +x^{*}_N(y_N), \qquad (y \in Y ), 
\]   
and let $a:E^N \times E^N \rightarrow \mathbb{R}$ be the continuous and bilinear form
\[
a(x,y):=a_1(x_1,y_1)+ \cdots + a_N(x_N,y_N), \qquad ((x,y) \in E^N \times E^N).
\]
First, let us verify that the problem \eqref{eq:1} is equivalent to finding an $x_0 \in \overline{Y}$ fulfilling the following condition: for each $y \in \displaystyle Y $, where $\overline{x}_0$ denotes the vector $(x_0,\dots,x_0)$,
\begin{equation}\label{eq:2}
x^{*}(y-\overline{x}_0)  \leq a(y,y-\overline{x}_0).
\end{equation}
The fact that (\ref{eq:1}) implies (\ref{eq:2}) follows from the sum of the inequalities and from the definition of $ x^{*}$ and $a$. For the opposite implication, it is enough to take $(y_1,x_0,\dots,x_0)$ as an element of $Y$ in (\ref{eq:2}) to obtain the first inequality of (\ref{eq:1}). We derive the other inequalities with the same reasoning.

Then we present the characterization mentioned at the beginning of this section, that extends the case of an equation previously established in \cite{pablo2020}:

\begin{theorem}\label{th:1}
Assume that $E$ is a real and reflexive Banach space, $Y_{1}, \dots, Y_{N}$ are nonempty and closed subsets of $E$, $x^{*}_1:E \longrightarrow \mathbb{R},  \dots, \, x^{*}_N:E \longrightarrow \mathbb{R}$ are continuous and linear functionals and define 
\[ 
x^{*}(y):= x^{*}_1(y_1)+ \cdots +x^{*}_N(y_N), \qquad \left( y \in E^N \right).
\] 
Let $a_1:E \times E \longrightarrow \mathbb{R}, \dots, \, a_N:E \times E \longrightarrow \mathbb{R}$ be continuous bilinear forms and le 
\[
a(x,y):=a_1(x_1,y_1)+ \cdots + a_N(x_N,y_N), \qquad ((x,y) \in E^N \times E^N).
\] 
Then, we have that there exists $\displaystyle x_0 \in \overline{Y}=\bigcap^N_{i=1} Y_i$ fulfilling the following system
\[
y \in Y \ \Rightarrow \ 
\
\left\{\begin{array}{lll}

x^{*}_1(y_1-x_0) & \leq & a_1(y_1,y_1-x_0)\\
x^{*}_2(y_2-x_0) & \leq & a_2(y_2,y_2-x_0)\\

\; \; \; \; \; \; \; \; \; \;  \vdots & \;  \vdots & \; \; \; \; \; \; \; \; \; \; \vdots\\

x^{*}_N(y_N-x_0) & \leq & a_N(y_N,y_N-x_0)\\

\end{array},
\right.
\]
if, and only if, for some $\alpha \geq 0$, $\overline{Y} \cap \alpha B_E \neq \emptyset$, and the next inequality holds:
\begin{equation}\label{eq:3}
\left.
	\begin{array}{c}
	m \ge 1, \ t \in \Delta_m  \\
	y_1, \dots,y_m \in Y
\end{array}	  
\right\}
\ \Rightarrow \
\sum_{j=1}^m t_j (x^*(y_j)-a(y_j,y_j))  \le \max_{x \in \overline{Y} \cap \alpha {B_E}} \left( \sum_{i=1}^N x_{i}^*(x)-a\left(\sum_{j=1}^m t_j y_j,\overline{x} \right) \right)  ,
\end{equation}
where $\overline{x}=(x,\dots,x)$.
\end{theorem}

\noindent \textsc{Proof.} We have that \eqref{eq:1} $\Rightarrow$ \eqref{eq:3} just by taking $\alpha:=\| x_0\|$.

For \eqref{eq:3} $\Rightarrow$ \eqref{eq:1}, let $X:= \overline{Y} \cap \alpha B_E$ and $Y=\displaystyle \prod_{i=1}^N Y_i$; choosing $m=1$, from \eqref{eq:3} we obtain
\[
\begin{array}{rl}
0 & \le \displaystyle \inf_{y \in Y} \left( a(y,y)-x^*(y) + \max_{x \in X} \left( \sum_{i=1}^N x_{i}^*(x)-a(y,\overline{x})\right) \right)    \\
  & = \displaystyle \inf_{y \in Y} \max_{x \in X} (a(y,y-\overline{x})-x^*(y- \overline{x})).
\end{array}
\]
Let
\[
\mu:= \displaystyle \inf_{y \in Y} \max_{x \in X} (a(y,y- \overline{x})-x^*(y- \overline{x})).
\]
If $\mu=-\infty$ there is nothing to prove. Otherwise, let $f:X \times Y \rightarrow \mathbb{R}$ be the function defined as
\[
f(x,y):=a(y,y- \overline{x})-x^{*}(y- \overline{x})-\mu, \; \; (\overline{x} \in X^N, y \in Y).
\]
From \eqref{eq:3} it follows
\[
\left.
	\begin{array}{c}
	m \ge 1, \ t \in \Delta_m  \\
	y_1, \dots,y_m \in Y
\end{array}	  
\right\}
\ \Rightarrow
0 \le \max_{x \in X} \sum_{j=1}^m t_j f(x,y_j),
\]
and as we have
\[
0=\displaystyle \inf_{y \in Y} \max_{x \in X} f(x,y),
\]
by Theorem \ref{th:minimax} the maximun is reached, and as a consequence, there exists $x_0 \in X$ satisfying
\[
\inf_{y \in Y} \max_{x \in X} f(x,y) \le \inf_{y \in Y} f(x_0,y),
\]
which is equivalent to the inequality system $\eqref{eq:1}$.
\hfill${\Box}$

\bigskip

Before stating our next result, we recall a technical lemma. Although it is proven \cite[Lemma 4.1]{cap14} for Hilbert spaces, it clearly works for Banach spaces.

\begin{lemma}\label{lm:1}
Let $E$ be a Banach space and let $Y$ be a nonempty closed convex subset of $E$. Assume that $a:E \times E \longrightarrow \mathbb{R}$ is a bilinear continuous form and $x^{\ast}:E \longrightarrow \mathbb{R}$ is a continuous and linear functional. Then, the next problem: find $y \in Y$ such that
\[
x^{*}(y-\overline{x}_0)  \leq a(y,y-\overline{x}_0),
\]
is equivalent to the problem of finding $y \in Y$ satisfying 
\[
x^{*}(y- \overline{x}_0)  \leq a(\overline{x}_0,y- \overline{x}_0).
\]
\end{lemma}

If we add certain more restrictive hypotheses to Theorem \ref{th:1}, we can equivalently express the condition \eqref{eq:3} in a simpler way, extending the system version of Stampacchia theorem.

\bigskip 

\begin{corollary}\label{co:1}
Let $E$ be a real and reflexive Banach space, let  $Y_{1}, \dots, Y_{N}$ be nonempty closed and convex sets of $E$, let $x^{*}_1:E \longrightarrow \mathbb{R},  \dots, \, x^{*}_N:E \longrightarrow \mathbb{R}$ be continuous and linear functionals and let $a_1:E \times E \longrightarrow \mathbb{R}, \dots, \, a_N:E \times E \longrightarrow \mathbb{R}$ be continuous bilinear forms. Let
\[ 
x^{*}(y):= x^{*}_1(y_1)+ \cdots +x^{*}_N(y_N), \qquad \left( y \in E^N \right),
\] 
and
\[
a(x,y):=a_1(x_1,y_1)+ \cdots + a_N(x_N,y_N), \qquad ((x,y) \in E^N \times E^N),
\] 
and suppose that
\[
x \in E \ \Rightarrow \ a(x,x) \geq0.
\]
Then, there exists $\displaystyle x_0 \in \overline{Y}=\bigcap^N_{i=1} Y_i$ such that
\begin{equation}\label{eq:4}
y \in Y:=\prod_{i=1}^N Y_i \; \Rightarrow \; \; x^{*}(y- \overline{x}_0)  \leq a(\overline{x}_0,y- \overline{x}_0)
\end{equation}
if, and only if, there exists $\alpha\geq0$ fulfilling $\overline{Y} \cap \alpha B_E \neq \emptyset$ and
\begin{equation}\label{eq:5}
y \in Y \; \Rightarrow \; \;  x^*(y)-a(y,y) \leq \displaystyle \sup_{x \in \overline{Y} \cap \alpha B_E} (x^*(\overline{x})-a(y,\overline{x})).
\end{equation}

\end{corollary}

\noindent \textsc{Proof}. First, let us observe that the conditions \eqref{eq:4} and \eqref{eq:2}, thanks to Lemma \ref{lm:1} with $u=(x_0,\dots,x_0)$ and $v=(y_1,\dots,y_N)$, are equivalents.  Therefore, we must prove that the conditions \eqref{eq:5} and \eqref{eq:3} are equivalents.

On the one hand, if we take $m=1$ in \eqref{eq:3} then we have \eqref{eq:5}. On the other hand, we suppose that \eqref{eq:5} is valid and let $m \geq 1$, $t \in \Delta_m$ and $y_1, \dots,y_m \in Y$. Then
\[
\begin{array}{rl}
	\displaystyle \sum_{j=1}^m t_j (x^*(y_j)-a(y_j,y_j)) & = \displaystyle  \sum_{j=1}^m t_j x^*(y_j) - \sum_{j=1}^m t_j a(y_j,y_j)   \\
	& \le \displaystyle x^* \left( \sum_{j=1}^m t_j y_j \right) -a \left(  \sum_{j=1}^m t_j y_j ,  \sum_{j=1}^m t_j y_j \right)   \\
	& \displaystyle \le \sup_{x \in \overline{Y} \cap \alpha {B_E}} \left( x^*(\overline{x})- a \left(  \sum_{j=1}^m t_j y_j , \overline{x} \right) \right),
\end{array}
\]
where we have used the convexity of $Y$ and that of the quadratic form associated with the bilinear form $a$ and \eqref{eq:5}.
\hfill${\Box}$

\bigskip

We conclude this section by proving that the version of systems of the classical Stampacchia theorem is a consequence of Corollary \ref{co:1}. Indeed, assume that $E$ is a real Hilbert space, $Y_1,\dots,Y_N$ are nonempty  closed and convex subsets of $E$, $x^{*}_1:E \longrightarrow \mathbb{R}, \dots, \, x^{*}_N:E \longrightarrow \mathbb{R}$ are continuous and linear functionals and $a_1:E \times E \longrightarrow \mathbb{R}, \dots, \,a_N:E \times E \longrightarrow \mathbb{R}$ are bilinear, continuous and coercive forms. With the notations above, let $x^{*}$ be the continuous and linear functional defined as
\[
 x^{*}(y):= x^{*}_1(y_1)+ \cdots +x^{*}_N(y_N), \qquad (y \in E^N),
\]   
and let $a:E^N \times E^N \rightarrow \mathbb{R}$ be the continuous and bilinear form
\[
a(x,y):=a_1(x_1,y_1)+ \cdots + a_N(x_N,y_N), \qquad ((x,y) \in E^N \times E^N).
\]
Let us note that, given a vector $x \in \alpha B_E$ with $\alpha \geq 0$, $\Vert x \Vert = \alpha$,  we can select, without loss of generality, the norm of $E^N$ appropriately so that $\Vert \overline{x} \Vert= \alpha$. 

In addition, if $\rho_1, \dots,\rho_N$ are the coercivity constants of $a_1, \dots, a_N$ and $x \in E$, then we have $(\rho_1+ \cdots+\rho_N) \Vert x \Vert^2 \leq a(x,x)$.

Let $\beta > 0$ such that $\overline{Y} \cap \beta B_E \neq \emptyset$ and $y \in Y:=\displaystyle \prod_{i=1}^N Y_i$. Then there hold
\[
\begin{array}{rl}
\displaystyle \frac{x^*(y)-a(y,y)}{\| y \|}-\frac{\displaystyle \sup_{x \in \overline{Y} \cap \beta {B_E}}(x^*(\overline{x})-a(y,\overline{x}))}{\| y \|} & \le \|x^*\|-(\rho_1+ \cdots+\rho_N) \| y \|+ \displaystyle \beta\frac{\|x^*-a(y,\cdot)\|}{\|y\|}   \\
 & \le \| x^* \| \left( 1 + \displaystyle \frac{\beta}{\| y \|}  \right)+ \beta \| a \|-(\rho_1+ \cdots+\rho_N) \| y \|. 
\end{array}
\]
Taking $\alpha > \beta$, it follows that $Y \cap \alpha B_E \neq \emptyset$ and
\[
y \in Y, \ \| y \| > \alpha \ \Rightarrow \ x^*(y)-a(y,y) \le \sup_{x \in \overline{Y} \cap \alpha {B_E}} (x^*( \overline{x} )-a(y, \overline{x}),
\]
and we arrive at \eqref{eq:5}.

\bigskip
\section{The inverse varational problem}

\noindent To solve the inverse problem associated with the system of variational inequalities \eqref{eq:1}, we will make use of the following collage-type result, which can be proved as a consequence of Stampacchia's theorem for a system of inequalities. In order to avoid expository complications, we previously introduce the following notation: for a real Banach space, $E^*$ is its topological dual space. Moreover, if $J$ is a nonempty set and for each $j \in J$ and $i \in \left\lbrace 1, \dots, N \right\rbrace $ ${x_j^1}^*,\dots,{x_j^N}^* \in E^*$ and $a_j^1,\dots ,a_j^N: E \times E \longrightarrow \mathbb{R}$ are continuos biliear forms, we denote by $({x_j^1}^*,\dots,{x_j^N}^*)$ and $(a_j^1,\dots,a_j^N)$ the continuous and linear functional
\[
{x_j^1}^*(y_1)+ \cdots +{x_j^N}^{*}(y_N), \qquad (y \in E^N),
\]   
and the continuous bilinear form
\[
a_j^1(x_1,y_1)+ \cdots + a_j^N(x_N,y_N), \qquad ((x,y) \in E^N \times E^N),
\]
respectively.

\begin{theorem}\label{th:2}
Let $J$ be a nonempty set, let $Y_1,\dots,Y_N$ be closed and convex nonempty subsebts of the Hilbert space $E$. For each $j \in J$ and $i \in \left\lbrace 1, \dots, N \right\rbrace $,  let ${x_j^1}^*,\dots,{x_j^N}^* \in E^*$ and $a_j^1,\dots,a_j^N: E \times E \longrightarrow \mathbb{R}$ be bilinear and continuous functionals satisfying that there exist $\rho^1_j,\dots,\rho^N_j$ positives such that
\[
y \in E \Rightarrow \ \rho^i_j \| y \|^2 \le a^i_j(y,y). 
\]
If $x_j^*:=({x_j^1}^*,\dots, {x_j^N}^*)$, $a_j:=(a^1_j,\dots,a^N_j)$, $Y:=Y_1 \times \cdots \times  Y_N$ and $\overline{x}_j$ is a solution of the system
\[
y \in Y \ \Rightarrow \ x_j^*(y-\overline{x}_j) \le a_j(\overline{x}_j,y-\overline{x}_j),
\]
then, 
\[
y \in Y, \ j \in J \  \Rightarrow \ \|y-\overline{x}_j \| \le \frac{\|a_j(y,\cdot)-x_j^*\|}{(\rho^1_j+\cdots+\rho^N_j)}.
\]
\end{theorem}

\noindent \textsc{Proof}. Given $y \in Y$ and $ j \in J$ and taking into account Corollary \ref{co:1}, we have
\[
\begin{array}{rl}
(\rho^1_j+\cdots+\rho^N_j) \| y-\overline{x}_j \|^2 & \le a_j(y-\overline{x}_j,y-\overline{x}_j)   \\
                   & = a_j(y,y-\overline{x}_j)-a_j(\overline{x}_j,y-\overline{x}_j)   \\
                   & \le a_j(y,y-\overline{x}_j)-x_j^*(y-\overline{x}_j)   \\
                   & = (a_j(y,\cdot)-x_j^*)(y-\overline{x}_j)   \\
                   & \le \| a_j(y,\cdot)-x_j^*\| \|y-\overline{x}_j\|.
\end{array}
\]
\hfill${\Box}$

In \cite{berkunlatrui,KuVr99,KuLaToVr06-sub,KuLaToLeRuGa15,KuLaMeVr} we can observe the idea that we used for the application of this result in the resolution of the inverse problem. This reasoning have been previously used with the Banach fixed point theorem in a similar way in \cite{kungom03}.

We finish our work with the following example, which consists of two clearly defined parts. The first deals with solving the forward problem using Galerkin's method. To this end, we will work with a certain Schauder basis, a very versatile tool, since we can observe its use in differential and integral problems (\cite{berkunlatrui, berforgarrui04, gamgarrui05,gamgarrui09, garrui2019,garrui2014}). The second part of this example, and also its main objective, is the numerical treatment of the inverse problem, where we obtain the target functions thanks to the Galerkin method previously described.

\begin{example}
Assume that $E:=H^1(0,1)$, $\lambda_1, \lambda_2$ are positive reals, $\alpha_1,\alpha_2,\beta_1,\beta_2 \in \mathbb{R}$ and $f,g \in L^\infty(0,1)$. We introduce the boundary value problem:

\[
\left\{\begin{array}{lll}

-u''(x)+ \lambda_1 u(x)=f(x) \quad \hbox{on } (0,1)\\

-v''(x)+ \lambda_2 v(x)=g(x) \quad \hbox{on } (0,1)\\

u(0)=\alpha_1, v(0)=\alpha_2, u(1)=\beta_1, v(1)=\beta_2

\end{array}
.
\right. 
\]
Considering the convex set
\[
Y:=\left\lbrace (w_1,w_2) \in E^2: w_1(0)=\alpha_1,w_2(0)=\alpha_2, w_1(1)=\beta_1, w_2(1)=\beta_2 \right\rbrace,
\]
and using a standard argumentation, we obtain the variational formulation of the previous system. Namely, for all $(w_1,w_2) \in Y$ it is satisfied that
\begin{equation}\label{eq:6}
\begin{array}{rl}
\displaystyle \int^1_0 u'(w_1-u)'+\int^1_0 v'(w_2-v)'+ \lambda_1 \int^1_0 u(w_1-u) + \lambda_2
\int^1_0 v(w_2-v) &  \geq \\
                   \displaystyle \int^1_0 f(w_1-u)+g(w_2-u).
\end{array}
\end{equation}
We define the bilineal, coercive and continuous form $a:E^2 \times E^2 \rightarrow \mathbb{R}$
\[
a((x_1,y_1),(x_2,y_2)):= \int^1_0 (x'_1 y'_1 +x'_2 y'_2) + \lambda_1 \int^1_0 x_1 y_1  + \lambda_2 \int^1_0 x_2 y_2, \; \; \; \; ((x_1,y_1),(x_2,y_2)) \in E^2 \times E^2,
\]
and the functional $x^{\ast}:E^2 \longrightarrow \mathbb{R}$
\[
x^{\ast}(x,y):=\int^1_0 (fx + gy) , \; \; ((x,y) \in E^2).
\]
The vectorial version of Stampacchia's theorem ensures the existence of a solution for the variational inequality \eqref{eq:6}. 

To show an example of the forward problem by using our numerical method, we define
\[
f(x):= \left( e - \dfrac{1}{5} \right) e^{\frac{x^2}{10}} - \frac{x^2}{25}e^{\frac{x^2}{10}}, \;\; (x \in \left[0,1\right] ),
\]
and
\[
g(x):=-2 \cos(x+1)^2+\frac{\pi}{2}\sin(x+1)^2+4(1+x)^2 \sin(x+1)^2 , \;\; (x \in \left[0,1\right] ).
\] 
We choose  $(\lambda_1,\lambda_2)=(e,\frac{\pi}{2})$. In order to use an appropriate Galerkin method, we take $z_1=w_1-u$ and $z_2=w_2-v$ in \eqref{eq:6}, and we obtain for each $(z_1,z_2) \in H^1_0 (0,1) \times H^1_0 (0,1)$
\[
\int^1_0 u'z'_1+\int^1_0 v'z'_2+ \lambda_1 \int^1_0 u z_1 + \lambda_2
\int^1_0 v z_2 \geq \int^1_0 fz_1+gz_2.
\]
To design the Galerkin method, we consider the Haar system  $\left\lbrace h_k  \right\rbrace_{k \geq 1} $ in $L^2(0,1)$. We define 
\[
g(x):=1, \;\; (x \in [0,1]),
\]
and for $k\geq2$
\[
g_k(x):=\int_0^x h_{k-1} (t) dt \;\; (x \in [0,1]).
\]
As a Schauder basis of $H^1(0,1)$ we use $\left\lbrace g_k  \right\rbrace_{k \geq 1}$, and for $H^1_0(0,1)$ we take $\left\lbrace g_{k+2}  \right\rbrace_{k \geq 1}$.

We have made use of Galerkin's method to solve the $m$-dimensional variational problem in the subspaces generated by $\left\lbrace g_3,\cdots,g_{m+2} \right\rbrace$. In the following table we show the behavior of the approximation in terms of the errors made in the corresponding spaces, where $(u_m, v_m)$ is the solution obtained for the $m$-dimensional problem. Also, we present some graphics that illustrate the exact solutions, their approximations and the differences between these functions.

\begin{figure}
 \centering
  \subfloat[Exact solution $u$ and their aproximations for $m=3,15$ and $63$]{
   \label{f:Envolventeconvexaconjunto}
    \includegraphics[width=0.5\textwidth]{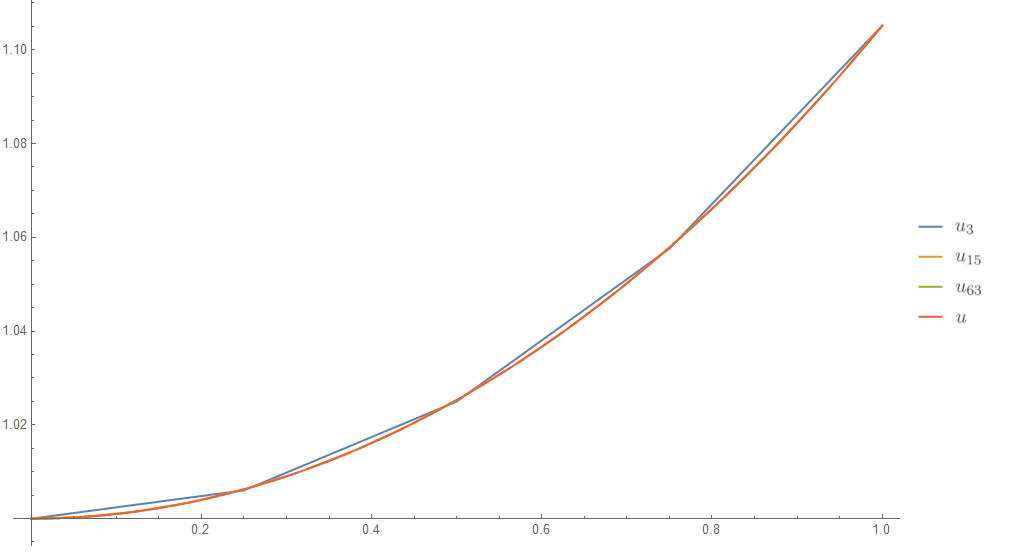}}
  \subfloat[Difference between $u$ and $u_{3}, u_{15}$, $u_{63}$]{
   \label{f:Conjuntodearco}
    \includegraphics[width=0.5\textwidth]{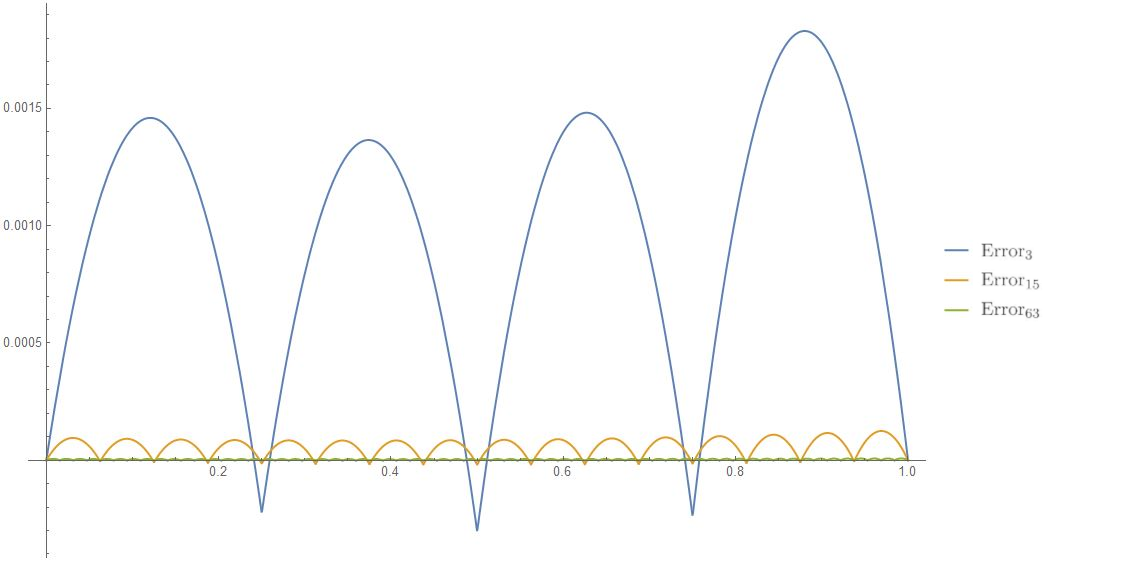}}
     \caption{}
 \label{f:Envolventeconvexadoble}
\end{figure}

\begin{figure}
 \centering
  \subfloat[Exact solution $u'$ and their aproximations for $m=3,15$ and $63$]{
   \label{f:Envolventeconvexaconjunto}
    \includegraphics[width=0.5\textwidth]{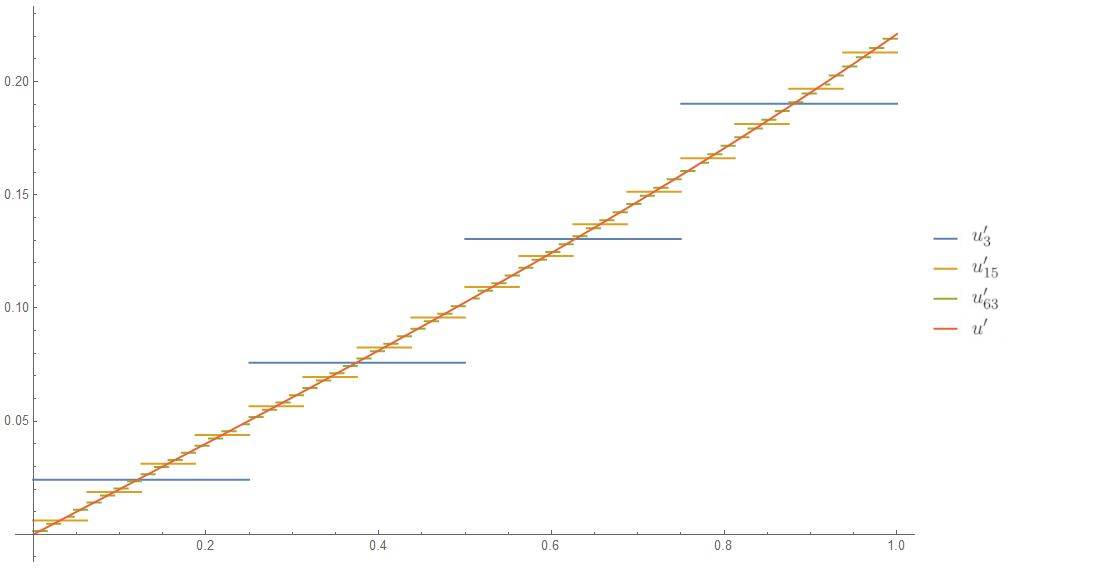}}
  \subfloat[Difference between $u'$ and $u'_{3}, u'_{15}$, $u'_{63}$]{
   \label{f:Conjuntodearco}
    \includegraphics[width=0.5\textwidth]{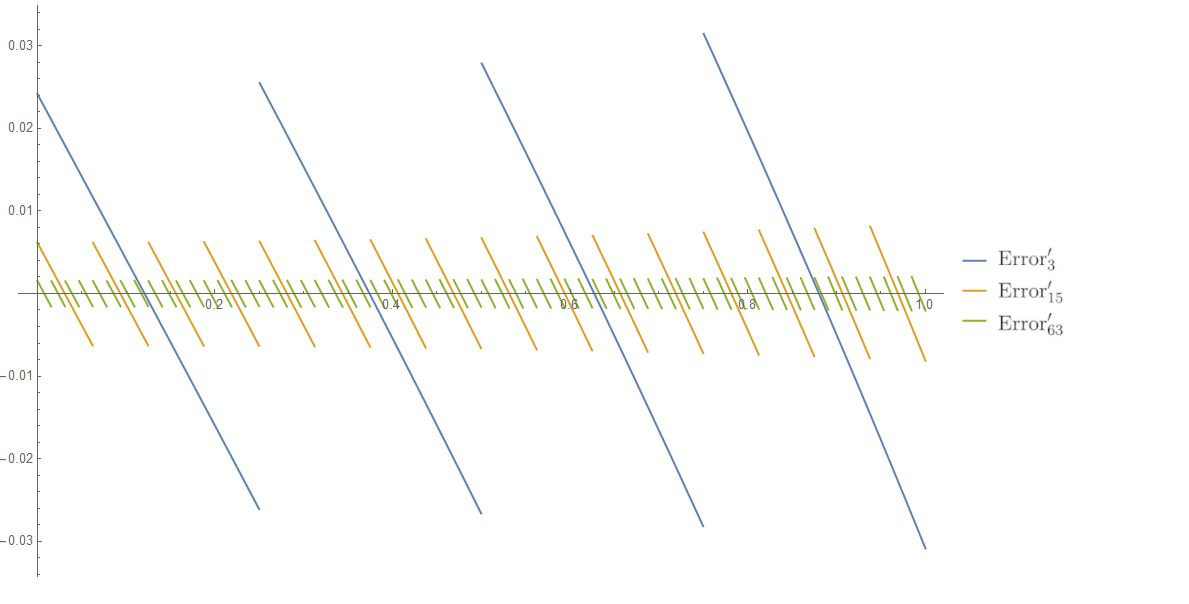}}
     \caption{}
 \label{f:Envolventeconvexadoble}
\end{figure}

\begin{figure}
 \centering
  \subfloat[Exact solution $v$ and their aproximations for $m=3,15$ and $63$]{
   \label{f:Envolventeconvexaconjunto}
    \includegraphics[width=0.5\textwidth]{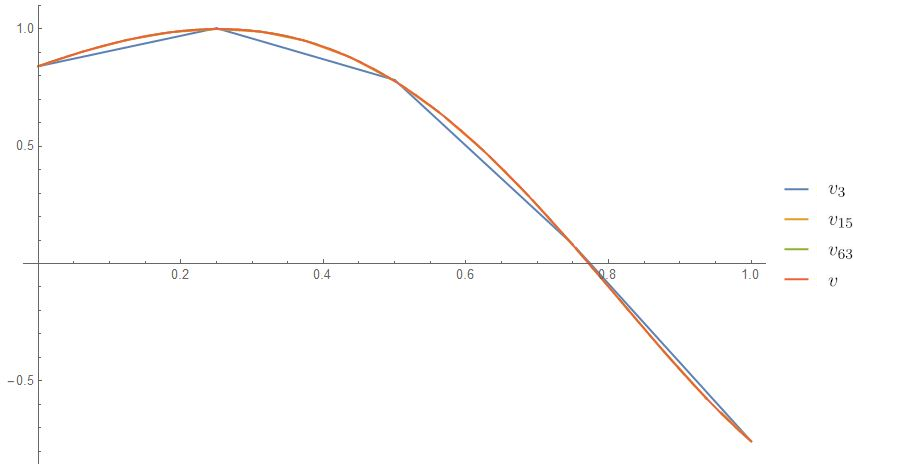}}
  \subfloat[Difference between $v$ and $v_{3}, v_{15}$, $v_{63}$]{
   \label{f:Conjuntodearco}
    \includegraphics[width=0.5\textwidth]{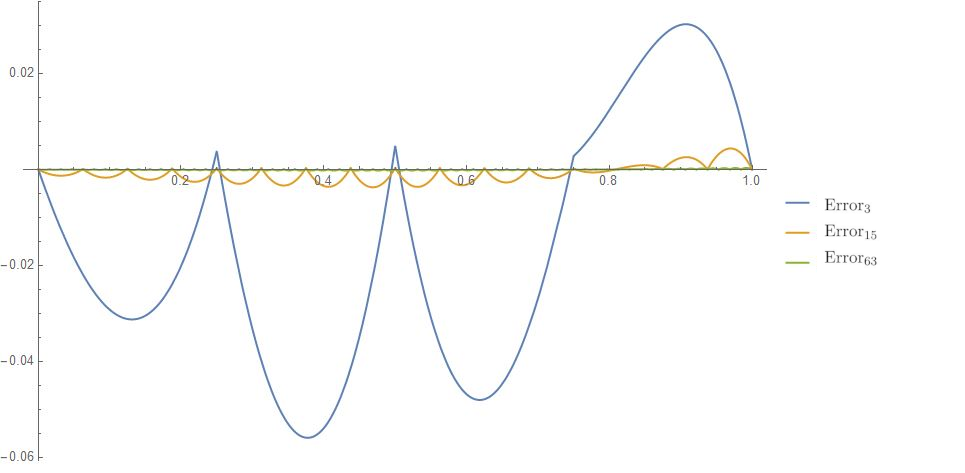}}
     \caption{}
 \label{f:Envolventeconvexadoble}
\end{figure}

\begin{figure}
 \centering
  \subfloat[Exact solution $v'$ and their aproximations for $m=3,15$ and $63$]{
   \label{f:Envolventeconvexaconjunto}
    \includegraphics[width=0.5\textwidth]{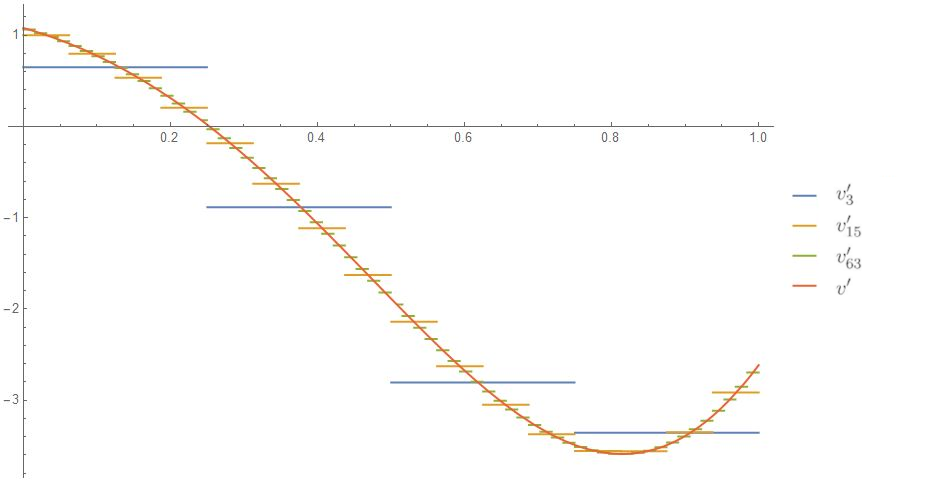}}
  \subfloat[Difference between $v'$ and $v'_{3}, v'_{15}$, $v'_{63}$]{
   \label{f:Conjuntodearco}
    \includegraphics[width=0.5\textwidth]{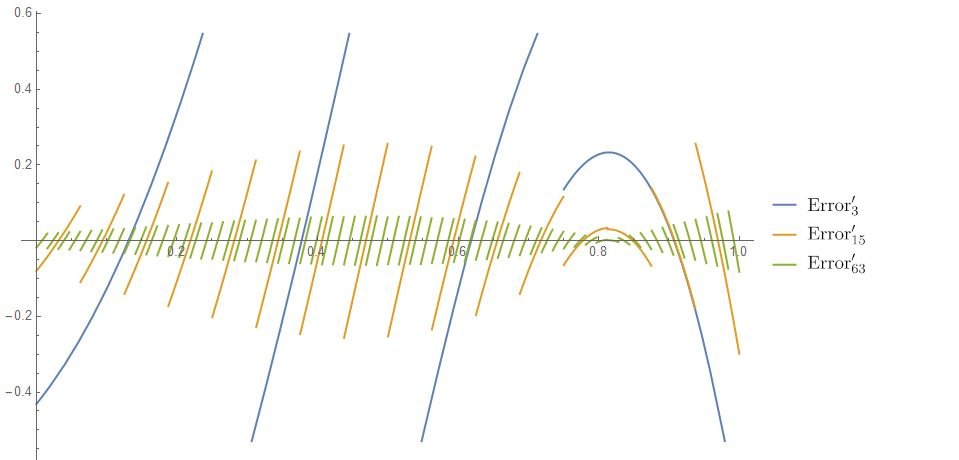}}
     \caption{}
 \label{f:Envolventeconvexadoble}
\end{figure}

\begin{table}[h]
\begin{center}
\begin{tabular}{|l|l|l|l|}
\hline
$m$ & $\Vert (u,v) - (u_{m},v_{m}) \Vert_{L_{2}}  $&   $\Vert (u',v') - (u'_m,v'_m) \Vert_{L_{2}} $ &  $ \Vert (u,v) - (u_m,v_m) \Vert_{H_{1}}  $ \\
\hline 
$3$ & $(0.00109781, 0.0304768) $ & $(0.0160122 , 0.421311 ) $ & $(0.0160498 , 0.422412 ) $ \\ \hline
$7$ & $(0.000272967, 0.00794639)$ & $(0.00800607 , 0.218486 ) $ & $ (0.00801072 , 0.21863) $ \\ \hline
$15$ & $(0.0000681497, 0.00200574 )$ & $( 0.004003, 0.110178) $  & $ (0.00400359 , 0.110196)$\\ \hline
$31$ & $(0.0000170317, 0.000502611 )$ & $ (0.0020015 , 0.0552043)$  & $(0.00200157 , 0.0552065 ) $\\ \hline
$63$ & $ (0.00000425756, 0.000125726 )$ & $(0.00100075 ,0.0276165 ) $  & $ (0.00100076 , 0.0276168)$\\ \hline
\end{tabular}
\caption{Errors of the forward method.}
\label{tabla:sencilla}
\end{center}
\end{table}

Now we finally present the treatment of the inverse problem under the notation of the Theorem \ref{th:2}. We now turn to the resolution of the inverse problem. The method we follow to solve it is as follows: we obtain a solution of the forward problem for each $j_0 \in J$. We write this solution as $x_j$. In the inverse problem, we try to find, whenever possible, a $j_0 \in J$ such that
\[
\| y-\overline{x}_{j_0}\|=\inf_{j \in J} \| y-\overline{x}_j\|.
\]
Let us take into account that if $Y$ is a closed affine subset of $E$, and $y \in Y$ is a target element, we can solve our problem if, under the condition that $\displaystyle \inf_{j}(\rho^1_j,\dots,\rho^N_j)>0$, we are able to solve
\begin{equation}\label{DonPablo}
\inf_{j \in J} \|(a_j(y,\cdot)-\overline{x}_j^*)_{| (Y-Y)}\|.
\end{equation}
To show a concrete example, we take the variational inequality discussed above, with $f(x):=\left( e - \dfrac{1}{5} \right) e^{\frac{x^2}{10}} - \frac{x^2}{25}e^{\frac{x^2}{10}}$, $g(x):=-2 \cos(x+1)^2+\frac{\pi}{2} \mathrm{sin}(x+1)^2+4(1+x)^2 \sin(x+1)^2$, $(\lambda_1,\lambda_2) \in ([0.5,3]\times[0.5,3])$ and $\alpha_1=1,\alpha_2=\sin 1,\beta_1=e^{\frac{1}{10}},\beta_2=\sin 4 $, that satisfy the hypotheses of Theorem 3.1 in a trivial way.

Now, taking $\lambda_1 = e$ and $\lambda_2 = \frac{\pi}{2}$ we obtain, using the forward method described previously, the approximate solutions $(u_m,v_m)$ for different values of $m$. We consider these approximate solutions as targets for solving the inverse problem. For our example, we can write \eqref{DonPablo} as
\[
\inf_{(\lambda_1,\lambda_2) \in ([0.5,3]\times[0.5,3])} \sup_{\scriptsize{\begin{array}{cc}\omega \in H^1_0(0,1)\\
\|\omega\|=1\end{array}}} |a_{j}((u_m,\omega),(v_m,\omega))-\overline{x}^*(\omega,\omega)|.
\]
We proceed to the discretization of our example, for that, we write the element $\omega \in H^1_0 (0,1)$ as a combination of the $n$ firsts elements of the Schauder basis of $H^1_0(0,1)$ given by $\left\lbrace g_{k+2}\right\rbrace_{k \geq 1}$ to obtain the minimization problem:
\[
 \inf_{(\lambda_1,\lambda_2) \in ([0.5,3]\times[0.5,3])}  \left|\sum_{k=1}^{n}(a_{j}((u_m,g_{k+2}),(v_m,g_{k+2}))-\overline{x}^*(g_{k+2},g_{k+2})\right|.
\]
Below we present a table where we can observe the different approximations of $(\lambda_1,\lambda_2)$ that we have obtained by taking $n=7$ in the above expression and considering different targets $(u_m,v_m)$.
\begin{table}[h]
\begin{center}
\begin{tabular}{|c|c|c|c|l|}
\hline
 $(u_3,v_3)$ &   $ (u_7,v_7) $ &  $ (u_{15},v_{15})  $ & $(u_{31},v_{31})$  \\
\hline  
 $ \tiny{(2.67488 ,2.95503)}$ & $ (2.70395 ,2.00997)$ & $(2.71828 ,1.5708)$ & $ (2.71828,1.5708)$   \\ \hline 
\end{tabular}
\caption{Numerical results for the inverse problem.}
\label{tabla:sencilla}
\end{center}
\end{table}

\end{example}

\bigskip

\section{Conclusions}

In this paper we have presented a numerical method to solve the inverse problem associated with a system of variational inequalities, which has been illustred in Example 3.2. To do this, firstly, we have used a minimax equality, Theorem 2.1, to prove a result that allows us to characterize the existence of a solution to the system of inequalities, Theorem 2.2. To solve the inverse problem, we have applied a consequence of Stampacchia's theorem, which is a collage-type result.

\bigskip

\section*{Acknowledgement}

Research partially supported by Junta de Andaluc\'{\i}a Grant FQM359.


\begin{thebibliography}{99}


\bibitem{ans99} Q. H. Ansari and J. C. Yao, \textsl{A fixed point theorem and its applications to a system of variational inequalities}, Bulletin of the Australian
Mathematical Society 59, (1999), 433--442.

\bibitem{aub1998} J.P. Aubin, \textsl{Optima and equilibria. An introduction to nonlinear analysis}, second edition, Graduate Texts in Mathematics 140, Springer-Verlag, Berlin, 1998.

\bibitem{berkunlatrui} M.I. Berenguer, H. Kunze, D. La Torre, M. Ruiz Gal\'an, \textsl{Galerkin method for constrained variational equations
and a collage-based approach to related inverse problems}, Journal of Computational and Applied Mathatematics 292 (2016), 67--75.

\bibitem{berforgarrui04} M.I. Berenguer, M.A. Fortes, A.I. Garralda-Guillem, M. Ruiz Gal\'an, \textsl{Linear Volterra integro-differential equation and Schauder bases}, Applied Mathematics and Computation 159 (2004), 495--507.

\bibitem{BiCasPa} G. Bigi, M. Castellani, M. Pappalardo, M. Passacantando, \textsl{Nonlinear Programming Techniques for Equilibria}, EURO advanced tutorials on operational research, Springer, Cham, 2019.

\bibitem{bofbrefor2013}  D. Boffi, F. Brezzi, M. Fortin, \textsl{Mixed finite element methods and applications}, Springer Series in Computational Mathematics 44, Springer, Heidelberg, 2013.




\bibitem{Capasso2014}
V. Capasso, H. Kunze, D. La Torre, E.R. and Vrscay, \textsl{Solving inverse problems for differential equations by a "generalized collage" method and
application to a mean field stochastic model}, Nonlinear Analysis: Real World Applications 15 (2014), 276--289.

\bibitem{cap14} A. Capatina, \textsl{Variational inequalities and frictional contact problems}, Advances in Mechanics and Mathematics 31, Springer, Cham, 2014.


\bibitem{cen12} Y. Censor, A. Gibali, S. Reich,  et al. \textsl{Common Solutions to Variational Inequalities}, Set-Valued Analysis 20, (2012), 229--247.

\bibitem{Fan1972} K. Fan, \textsl{A minimax inequality and applications}, Inequalities, III (Proc. Third Sympos. Univ. California, Los Angeles, CA, 1969; dedicated to the memory of Theodore S. Motzkin), 103--113, Academic Press, New York, 1972.

\bibitem{fan} K. Fan, \textsl{Minimax theorems},
Proceedings of the National Academy of Sciences of the United States of America \textbf{39} (1953), 42--47.

\bibitem{gamgarrui05} D. G\'amez, A.I. Garralda-Guillem, M. Ruiz Gal\'an, \textsl{Nonlinear initial-value problems and Schauder bases}, Nonlinear Analysis 63 (2005), 97--105.

\bibitem{gamgarrui09} D. G\'amez, A.I. Garralda-Guillem, M. Ruiz Gal\'an, \textsl{High-order nonlinear initial-value problems countably determined}, Journal of Computational and Applied Mathatematics 228 (2009), 77--82.


\bibitem{garrui2019} A.I. Garralda-Guillem, M. Ruiz Galán, \textsl{A minimax approach for the study of systems of variational equations and related Galerkin schemes}, Journal of Computational and Applied Mathatematics 354 (2019), 103--111.


\bibitem{garrui2014} A.I. Garralda-Guillem, M. Ruiz Galán, \textsl{Mixed variational formulations in locally convex spaces}, Journal of Mathematical Analysis and Applications 414 (2014), 825--849.



\bibitem{gat2014} G.N. Gatica, \textsl{A simple introduction to the mixed finite element method},  Theory and applications, SpringerBriefs in Mathematics, Springer, Cham, 2014.


\bibitem{kas00} G. Kassay and J. Kolumbán, \textsl{System of multi-valued variational inequalities}, Publicationes Mathematicae Debrecen 56 (2000), 185--195.

\bibitem{kaskol1996} G. Kassay and J. Kolumb\'an, \textsl{On a generalized sup-inf problem}, Journal of Optimization Theory and Applications \textbf{91} (1996), 651--670.



\bibitem{kon97} I. V. Konnov, \textsl{On systems of variational inequalities}, Russian Mathematics 41 (1997), 79--88.








\bibitem{KuVr99}
H. Kunze, E.R. Vrscay, \textsl{Solving inverse problems for ordinary differential equations using the Picard contraction mapping}, Inverse Problems 15 (1999), 745--770.

\bibitem{kunze03a}
H. Kunze, S. Gomes, \textsl{Solving An Inverse Problem for Urison-type Integral Equations Using Banach's Fixed Point Theorem}, Inverse Problems 19 (2003), 411--418.

\bibitem{kunze03b}
H. Kunze, J. Hicken, E.R. Vrscay, \textsl{Inverse Problems for ODEs Using Contraction Maps: Suboptimality of the "Collage Method"}, Inverse Problems 20 (2004), 977--991.


\bibitem{KuLaToVr06-sub}
H. Kunze, D. La Torre, E. R. Vrscay,\textsl{ A generalized collage method based upon the Lax--Milgram functional for solving boundary value inverse problems},  Nonlinear Analysis 71 (2009), e1337--e1343.

\bibitem{Ku1}
H. Kunze, D. La Torre, E. R. Vrscay, \textsl{Solving inverse problems for variational equations using the "generalized collage methods," with applications to boundary value problems}, Nonlinear Analysis Real World Applications, 11 (2010), 3734--3743.



\bibitem{KuLaToLeRuGa15}
H. Kunze, D. La Torre, K. Levere, M. Ruiz Galan, \textsl{Inverse problems via the "generalized collage theorem" for vector-valued Lax-Milgram-based variational problems}, Mathematical Problems in Engineering, (2015), 1--8.

\bibitem{KuLaMeVr}
H. Kunze, D. La Torre, F.Mendivil, E. R. Vrscay, \textsl{Fractal-based methods in analysis}, Springer, 2012.



\bibitem{kungom03} H. Kunze, G. Gomes, \textsl{Solving an inverse problem for Urison-type integral equations using Banach's fixed point theorem}, Inverse Problems 19 (2003), 411--418.

\bibitem{Levere}
K. Levere, H. Kunze, D. La Torre, \textsl{A collage-based approach to solving inverse problems for second-order nonlinear parabolic PDEs}, Journal of Mathematical Analysis and Applications, 406 (2013), 120--133.

\bibitem{pablo2020} P. Montiel López, \textsl{A minimax approach for inverse variational inequalities}, Communications in Nonlinear Science and Numerical Simulation 90, Article number 105339.



\bibitem{pan}
J. S. Pang, \textsl{ Asymmetric variational inequality problems over product sets: Applications and iterative methods}, Mathematical Programming 31, (1985), 206--219.

\bibitem{Par2018}
S. Park, \textsl{On the minimax inequality of Brezis-Nirenberg-Stampacchia}, Journal of Nonlinear and Convex Analysis 19(9), (1985), 1493--1501.


\bibitem{rui2014} M. Ruiz Galán, \textsl{An intrinsic notion of convexity for minimax}, Journal of Convex Analysis 21 (2014), 1105--1139. 

\bibitem{rui2016} M. Ruiz Galán, \textsl{The Gordan theorem and its implications for minimax theory}, Journal of Nonlinear and Convex Analysis 17 (2016), 2385--2405. 


\bibitem{sai2018} J. Saint Raymond, \textsl{A new minimax theorem for linear operators}, Minimax Theory and its Applications 3 (2018), 131--160. 

\bibitem{sim1998} S. Simons, \textsl{Minimax and monotonicity}, Lecture Notes in Mathematics 1693, Springer, New York, 1998.

\bibitem{sim2007} S. Simons, \textsl{The Hahn--Banach--Lagrange theorem},  Optimization 56 (2007), 149--169.

\bibitem{ste} A. Stefanescu, \textsl{A theorem of the alternative and a two-function minimax
theorem}, Journal of Applied Mathematics \textbf{2004:2} (2004),
167--177.





\bibitem{zao10} Y. Zhao, Z. Xia, L. Pang and L. Zhang, \textsl{Existence of solutions and
algorithm for a system of variational inequalities}, Fixed Point Theory
and Applications (2010), Article ID 182539, 11 pp.

\end{thebibliography}
\end{document}